\documentclass[12pt]{amsart}
\usepackage{graphicx}
\usepackage{amssymb,amsmath,amscd,amsthm,pb-diagram}

\newcommand{\Z}{{\mathbb Z}}

\DeclareMathOperator{\fr}{fr}

\DeclareMathOperator{\pr}{pr}

\newtheorem*{maintheorem}{Theorem}
\newtheorem*{question}{Question}
\newtheorem*{mainlemma}{Lemma}

\theoremstyle{definition}

\newtheorem*{acknowledgment}{Acknowledgment}

\title[Concordance to the Hopf link]{A two component link with Alexander polynomial one is concordant to the Hopf link}
\author{James F. Davis}
\thanks{Partially supported by the NSF}
\address{ Department of Mathematics\\
Indiana University\\
Bloomington, IN 47405\\
USA}
\email{jfdavis@indiana.edu}
\date{}

\begin{document}

\begin{abstract}   Topological surgery in dimension four is used to show a two component link with Alexander polynomial one is topologically concordant to the Hopf link.
\end{abstract}
\maketitle

Let $L$ be a two component link in $S^3$, an embedding of two disjoint circles which is topologically locally flat, that is, which extends to an embedding of two solid tori.  The link has Alexander polynomial one, if, and only if, the first homology of the universal abelian cover of the complement of the link vanishes.  If the link has Alexander polynomial one, then the linking number of the two components is one.

Two links are {\em concordant} if there is an topologically locally flat embedding
$$
(S^1 \amalg S^1) \times I \to S^3 \times I
$$
which restricts to the given links 
$$
(S^1 \amalg S^1) \times \{i\} \to S^3 \times \{i\}
$$
when $i = 0,1$.

\begin{maintheorem}  A two component link with Alexander polynomial one is concordant to the Hopf link. 
\end{maintheorem}

Freedman showed that a knot with Alexander polynomial one is concordant to the unknot \cite{F}, \cite[p. 210]{FQ}.  Hillman had a program  for proving the above theorem (see \cite[Section 7.6]{H2}, which corrected the account in \cite{H1}), but asked whether a final obstruction must vanish.  In this note we complete Hillman's program and show that the surgery problem can be chosen so that the final obstruction vanishes. Topological surgery in dimension four is used; note that the fundamental group of the complement of the Hopf link is $\Z^2$, which is a good group in the sense of \cite{FQ}.  

Casson (unpublished) observed that Donaldson's results on gauge theory lead to examples of Alexander polynomial one knots which are not smoothly concordant to the unknot; Fintushel-Stern \cite{FS} use gauge theory to show that the $(-3,5,7)$ pretzel knot is such a knot.
It follows that our theorem is not true in the smooth category.  Here is a more refined question.

\begin{question}
Is there a two component link with Alexander polynomial one which is not smoothly concordant to the Hopf link, but each of whose components are smoothly concordant to the unknot?
\end{question}

\section{Proof of Theorem}
The strategy of the proof is to first build a closed 3-manifold $M_L$, which will eventually serve as the boundary of the complement of a tubular neighborhood of the embedding of $(S^1 \amalg S^1) \times I$.  We then show this 3-manifold is a framed boundary.  We surger the 4-manifold and then glue in $(S^1 \amalg S^1) \times I \times D^2$ to arrive at the desired concordance.

Let $H$ be the Hopf link and $L$ be a two component link with Alexander polynomial one.  Let
$$
M_L = (S^3 - N(L)) \cup_{T^2\amalg T^2} ((T^2 \amalg T^2) \times I) \cup_{T^2\amalg T^2} (S^3 - N(H)), 
$$
where $N(L)$ and $N(H)$ are tubular neighborhoods of the links, and $(T^2 \amalg T^2) \times 0$ and $(T^2 \amalg T^2) \times 1$ are identified with the boundary of the closures of $N(L)$ and $N(H)$ respectively, so that one factor of each $T^2 = S^1 \times S^1$ maps to the meridian and the other factor maps to a longitude whose linking number with the corresponding component of the link is zero.  Then $M_L$ is a framed 3-manifold and $H_1(M_L) = \Z^3$.  Note that $M_H$ is the 3-torus.
%o\footnote{Comment this out later.  By collapsing out the $(T^2 \amalg T^2) \times I$, one sees that $M_L$ is homeomorphic to the union of $S^3 - N(L)$ and $S^3 - N(H)$ along their common boundary.  Also, in the case where $L = H$, $M_L = T^3$.}

Choose a map $c : M_L \to T^2$ which restricts to homology isomorphisms on $S^3 - N(L)$ and $S^3 - N(H)$.  

\begin{mainlemma}  There is a framing on $M_L$ so that $[c : M_L \to T^2] \in \Omega_3^{\fr}(T^2)$ represents the zero element.
\end{mainlemma}

\begin{proof}  There is a geometrically defined map 
\begin{align*}
\Omega^{\fr}_3(T^2) \to & \Omega^{\fr}_3 \times \Omega^{\fr}_2 \times \Omega^{\fr}_2 \times \Omega^{\fr}_1 = \Z_{24} \times \Z_2 \times \Z_2 \times \Z_2\\
[f : M \to T^2] \to &([M], [(\pr_1 \circ f)^{-1}*], [(\pr_2 \circ f)^{-1}*],[ f^{-1}*]),
\end{align*}
where $\pr_i : T^2 \to S^1$ are projection onto the factors, and all maps are assumed to be transverse to $*  \in S^1$ or $* \in T^2$.  One can see this map is an isomorphism either by the collapsing of the Atiyah-Hirzebruch spectral sequence or by seeing geometrically that the map is onto and by using a cofibration argument to show that domain has at most 192 elements.

 The J-homomorphism $\pi_3 O \to \pi_3^S = \Omega^{\fr}_3$ is onto, so one can alter the framing of $M_L$ in a disk neighborhood of a point to change  the first component $[c : M_L \to T^2]$ to zero.  Let $N^1 = f^{-1}*$.  Change the framing on a normal bundle of $N^1$ by a map to $SO(2)$ and extend this map to $M_L$ using that $H^1(M_L) \to H^1(N^1)$ is surjective.
 %Use that $M_L$ is a homology 3-torus.
Thus one can assume the last component of $[c]$ is trivial too.

The image of $[c]$ in $ \Omega^{\fr}_2 \times \Omega^{\fr}_2$ is given by the Arf invariant of each component of the link, and these vanish because each component has Alexander polynomial one, and hence Arf invariant one, by Levine's criterion \cite {L} that the Arf invariant of a knot vanishes if and only if $\Delta_K(-1) \equiv \pm 1 \pmod 8$. (Alternatively, Freedman-Quinn \cite[p. 207] {FQ} give a direct geometric argument for the vanishing of the Arf invariant for a knot of polynomial one.)  \end{proof}

Let $W^4_L$ be a compact framed 4-manifold with boundary $M_L$ and $f: W^4_L \to T^2$ a map which restricts on the boundary to $c$.  After performing surgeries on  circles generating the kernel of $\pi_1(f)$, we may assume that $f$ induces an isomorphism on $\pi_1$.  By using obstruction theory,  find a degree one normal map
$$F : (W_L, M_L) \to (T^2 \times I \times I, T^3)$$
so that $\pr_{T^2} \circ F = f : W_L \to T^2$ and so that $F$ restricted to the boundary is a $\Z[\Z^2]$-homology equivalence.   There is an obstruction $\theta(F) \in L_4(\Z[\Z^2]) = L_4(\Z) \oplus L_2(\Z) = \Z \oplus \Z_2$   to doing surgery rel $\partial$ to a map $G : (X, M_L) \to  (T^2 \times I \times I, T^3)$ so that $G :X \to T^2 \times I \times I$ is a homotopy equivalence.  The surgery obstruction $\theta(F)$ is detected geometrically by the signature of $W_L$ divided by 8, and the Arf invariant of $F^{-1}(* \times I^2,* \times \partial I^2) \to (* \times I^2,* \times \partial I^2)$, after perturbing $F$ so that it is transverse to $(* \times I^2,* \times \partial I^2)$ and so that $F^{-1}(* \times \partial I^2) \to (* \times \partial I^2)$ is a homotopy equivalence (see \cite{S}).  This can be arranged in our case.  We can then modify $F$ and $W_L$ by first connecting sum with $\sigma(W_L)/8$-copies of Freedman's $E_8$-manifold to kill the signature invariant. If there is a nontrivial Arf invariant,  find a closed surface $\Sigma^2$ in $W_L$ representing the generator of $H_2(T^2 \times I \times I)$ and replace its neighborhood $\Sigma^2 \times D^2$ with $\Sigma^2 \times (T^2 - \text{int }D^2)$ where $T^2$ is given the Arf invariant framing.  This  alters\footnote{This alteration of the surgery obstruction can also be seen algebraically by identifying the normal invariants relative to the boundary with $[(T^2 \times I^2, T^2 \times \partial I^2), (G/TOP,*)] = H^4(T^2 \times I^2, T^2 \times \partial I^2;\Z) \oplus H^2(T^2 \times I^2, T^2 \times \partial I^2;\Z_2)= \Z \oplus \Z_2$. This maps isomorphically to $L_4(\Z[\Z^2])$ via the composition of duality with the assembly map.} the surgery problem $F$ to one where surgery is possible and produces a map $G : (X, M_L) \to  (T^2 \times I \times I, T^3)$ so that $G :X \to T^2 \times I \times I$ is a homotopy equivalence.  

By gluing in a product of the interval with a  neighborhood of the link, one arrives at the desired concordance $$X \cup_{(T^2 \amalg T^2) \times I} (S^1 \times D^3 ~ \amalg ~ S^1 \times D^3) \cong S^3 \times I,$$ where the homeomorphism is due to the proof of the four-dimensional Poincar\'e conjecture.  \qed
\\

\noindent {\bf Dialogue} The following conversation was overheard:
\\

\noindent {\em Low-dimensional topologist}:  The exposition of the surgery theory in the argument above seems somewhat terse.  For example, where was the Alexander polynomial one hypothesis used?
\\

\noindent {\em Surgeon}:  Hmm.  To make sense of the rel $\partial$ surgery obstruction $\theta(F) \in L_4(\Z[\Z^2])$, one needs that $F|_{M_L} : M_L \to T^3$ is a $\Z[\pi_1(T^2 \times I \times I)]$-homology equivalence.  The vanishing of the first homology of the universal abelian cover of the link complement implies that the $\Z[\Z^2]$-cover of $M_L$ is a homology circle.
\\

\noindent{\em Low-dimensional topologist}:  What's new in this paper?  For example, how does it compare with Hillman's attack's on the problem in \cite{H2}?
\\

\noindent{\em Surgeon}:  Well, Jonathan constructed a similar surgery problem with obstruction in $L_4(\Z[\Z^2]) = \Z \oplus \Z_2$, but did not know if the $\Z_2$-component was zero.  The main contribution here is the realization that the surgery problem could be altered.

\begin{acknowledgment}  The part of Low-dimensional topologist was played by Charles Livingston, with understudy Stefan Friedl.
\end{acknowledgment}

 \end{document}